# SECTIONAL GENERA OF PARAMETER IDEALS

SHIRO GOTO AND KAZUHO OZEKI

Abstract. Let $M$ be a finitely generated module over a Noetherian local ring. This paper reports, for a given parameter ideal $Q$ for $M$, a criterion for the equality $\mathrm{g}_s(Q;M) = \mathrm{hdeg}_Q(M) - \mathrm{e}_Q^0(M) - \mathrm{T}_Q^1(M)$, where $\mathrm{g}_s(Q;M)$, $\mathrm{e}_Q^0(M)$, $\mathrm{e}_Q^1(M)$, and $\mathrm{T}_Q^1(M)$ respectively denote the sectional genus, the multiplicity, the first Hilbert coefficient, and the Homological torsion of $M$ with respect to $Q$.

Contents

## 1. Introduction

The notion of the sectional genera of commutative rings was introduced by A. Ooishi [O], and since then, many authors have been engaged in the development of the theory. The purpose of our paper is to give a criterion for a certain equality of the sectional genera of parameters for modules.

To state the problems and the results of our paper, let us fix some of our terminology. Let $A$ be a Noetherian local ring with maximal ideal $\mathfrak{m}$ and $d = \dim A > 0$. Let $M$ be a finitely generated $A$-module with $s = \dim_A M$. For simplicity, throughout this paper, we assume that $A$ is $\mathfrak{m}$–adically complete and the residue class field $A/\mathfrak{m}$ of $A$ is infinite. Let $I$ be a fixed $\mathfrak{m}$-primary ideal in $A$ and let $\ell_A(N)$ denote, for an $A$-module $N$, the length of $N$. Then there exist integers $\{\mathrm{e}_I^i(M)\}_{0 \le i \le s}$ such that

$$\ell_A(M/I^{n+1}M) = \mathrm{e}_I^0(M)\binom{n+s}{s} - \mathrm{e}_I^1(M)\binom{n+s-1}{s-1} + \cdots + (-1)^s\mathrm{e}_I^s(M)$$

for all $n \gg 0$. We call $\mathrm{e}_I^i(M)$ the $i$-th Hilbert coefficient of $M$ with respect to $I$ and especially call the leading coefficient $\mathrm{e}_I^0(M)$ $(> 0)$ the multiplicity of $M$ with respect to $I$. We set

$$\mathrm{g}_s(I;M) = \ell_A(M/IM) - \mathrm{e}_I^0(M) + \mathrm{e}_I^1(M)$$

and call it the sectional genus of $M$ with respect to $I$.

*Key words and phrases:* Hilbert function, Hilbert coefficient, homological degrees

In this paper we need the notions of homological degrees and torsions of modules. For each $j \in \mathbb{Z}$ we set

$$M_j = \mathrm{Hom}_A(\mathrm{H}^j_{\mathfrak{m}}(M), E),$$

where $E = \mathrm{E}_A(A/\mathfrak{m})$ denotes the injective envelope of $A/\mathfrak{m}$ and $\mathrm{H}^j_{\mathfrak{m}}(M)$ the $j$th local cohomology module of $M$ with respect to the maximal ideal $\mathfrak{m}$. Then $M_j$ is a finitely generated $A$-module with $\dim_A M_j \leq j$ for all $j \in \mathbb{Z}$ (Fact 2.1).

The homological degree $\mathrm{hdeg}_I(M)$ of $M$ with respect to $I$ is inductively defined in the following way, according to the dimension $s = \dim_A M$ of $M$.

**Definition 1.1.** ([V2]) For each finitely generated $A$-module $M$ with $s = \dim_A M$, we set

$$\mathrm{hdeg}_I(M) = \begin{cases} \ell_A(M) & \text{if } s \leq 0, \\ & \\ \mathrm{e}^0_I(M) + \sum_{j=0}^{s-1} \binom{s-1}{j} \mathrm{hdeg}_I(M_j) & \text{if } s > 0 \end{cases}$$

and call it the homological degree of $M$ with respect to $I$.

The homological torsion of $M$ with respect to $I$ is defined as follows.

**Definition 1.2.** Let $M$ be a finitely generated $A$–module with $s = \dim_A M \geq 2$. We set

$$\mathrm{T}^i_I(M) = \sum_{j=1}^{s-i} \binom{s-i-1}{j-1} \mathrm{hdeg}_I(M_j)$$

for each $1 \leq i \leq s-1$ and call them the homological torsions of $M$ with respect to $I$.

Notice that the homological degrees $\mathrm{hdeg}_I(M)$ and torsions $\mathrm{T}^i_I(M)$ of $M$ with respect to $I$ depend only on the integral closure of $I$.

Let $Q = (a_1, a_2, \ldots, a_s)$ be a parameter ideal for $M$. We denote by $\mathrm{H}_i(Q; M)$ $(i \in \mathbb{Z})$ the $i$–th homology module of the Koszul complex $\mathrm{K}_\bullet(Q; M)$ generated by the system $a_1, a_2, \ldots, a_s$ of parameters for $M$. We set

$$\chi_1(Q; M) = \sum_{i \geq 1} (-1)^{i-1} \ell_A(\mathrm{H}_i(Q; M))$$

and call it the first *Euler characteristic* of $M$ relative to $Q$; hence

$$\chi_1(Q; M) = \ell_A(M/QM) - \mathrm{e}^0_Q(M) \geq 0$$

by a classical result of Serre (see [AB], [Se]).

In [GhGHOPV, Theorem 7.1], it was proved that, for parameter ideals $Q$ for $M$, the upper bound $\chi_1(Q; M) \leq \mathrm{hdeg}_Q(M) - \mathrm{e}^0_Q(M)$ of the first Euler characteristic $\chi_1(Q; M)$ of $M$ relative to $Q$. In [GO, Theorem 1.3], the authors gave a criterion for the equality $\chi_1(Q; M) = \mathrm{hdeg}_Q(M) - \mathrm{e}^0_Q(M)$. We also have the inequalities $0 \geq \mathrm{e}^1_Q(M) \geq -\mathrm{T}^1_Q(M)$ for every parameter ideals $Q$ for $M$ ([MSV, Theorem 3.6], [GhGHOPV, Theorem 6.6]), where the equality $\mathrm{e}^1_Q(M) = -\mathrm{T}^1_Q(M)$ holds true if and only if the equality $\chi_1(Q; M) = \mathrm{hdeg}_Q(M) - \mathrm{e}^0_Q(M)$ holds true ([GO, Theorem 1.4]), provided $M$ is unmixed. The reader may consult [GhGHOPV] for the characterization of modules which have parameter ideals $Q$ with $\mathrm{e}_Q^1(M) = 0$. Thus the behavior of the first Euler characteristics $\chi_1(Q;M)$ and the first Hilbert coefficients $\mathrm{e}_Q^1(M)$ are rather satisfactory understood.

In this paper we study the sectional genus

$$\mathrm{g}_s(Q;M) = \ell_A(M/QM) - \mathrm{e}_Q^0(M) + \mathrm{e}_Q^1(M)$$

of $M$ with respect $Q$ in connection with homological degrees and torsions. First, we will show that in the case where $\dim_A M = 1$ the inequality

$$\mathrm{g}_s(Q;M) \leq 0$$

holds true for every parameter ideals $Q$ for $M$ (Lemma 3.1). We will also show that $\mathrm{g}_s(Q;M) = 0$ if and only if the ideal $Q$ is generated by a parameter $a$ for $M$ which forms a $d$-sequence on $M$ (Lemma 3.1). We should note that, in [GHV, Mc], it was proved that the inequality $\chi_1(Q;A) \leq -\mathrm{e}_Q^1(A)$ holds true for parameter ideals $Q$ in a Noetherian local ring $A$ with depth$A \geq d-1$. In [GHV], Hong, Vasconcelos, and the first author gave a criterion for the equality $\chi_1(Q;A) = -\mathrm{e}_Q^1(A)$.

Let $M$ be a finitely generated $A$-module with $\dim_A M \geq 2$. Then we have the inequality

$$\mathrm{g}_s(Q;M) \leq \mathrm{hdeg}_Q(M) - \mathrm{e}_Q^0(M) - \mathrm{T}_Q^1(M) = \sum_{i=0}^{s-2} \binom{s-2}{i} \mathrm{hdeg}_Q(M_i)$$

for all parameter ideals $Q$ for $M$ (Proposition 3.3). Hence the upper bound of $\mathrm{g}_s(Q;M) = \chi_1(Q;M) + \mathrm{e}_Q^1(M)$ is given by the sum of the upper bound of $\chi_1(Q;M)$ and the lower bound of $\mathrm{e}_Q^1(M)$. It seems natural to ask what happens on the parameters $Q$ for $M$, when the equality

$$\mathrm{g}_s(Q;M) = \mathrm{hdeg}_Q(M) - \mathrm{e}_Q^0(M) - \mathrm{T}_Q^1(M)$$

is attained.

The main result of this paper is stated as follows, where the sequence $a_1, a_2, \ldots, a_d$ is said to be a $d$-sequence on $M$, if the equality

$$[(a_1, a_2, \ldots, a_{i-1})M :_M a_i a_j] = [(a_1, a_2, \ldots, a_{i-1})M :_M a_j]$$

holds true for all $1 \leq i \leq j \leq d$ ([H]).

**Theorem 1.3.** *Let $M$ be a finitely generated $A$-module with $d = \dim_A M \geq 2$ and let $Q$ be a parameter ideal of $A$. Then the following conditions are equivalent:*

(1) $\mathrm{g}_s(Q;M) = \mathrm{hdeg}_Q(M) - \mathrm{e}_Q^0(M) - \mathrm{T}_Q^1(M)$.

(2) *The following two conditions are satisfied:*

(*a*)

$$(-1)^i \mathrm{e}_Q^i(M) = \begin{cases} \mathrm{T}_Q^i(M) & \text{if } 2 \leq i \leq d-1, \\ \ell_A(\mathrm{H}_\mathfrak{m}^0(M)) & \text{if } i = d \end{cases}$$

*for all $2 \leq i \leq d$.*

(*b*) $\ell_A(M/QM) - \sum_{i=0}^{d} (-1)^i \mathrm{e}_Q^i(M) = 0$.

*When this is the case, we have the following*:

(i) *there exist elements* $a_1, a_2, \ldots, a_d \in A$ *such that* $Q = (a_1, a_2, \ldots, a_d)$ *and* $a_1, a_2, \ldots, a_d$ *forms a d-sequence on* $M$,

(ii)
$$\ell_A(M/Q^{n+1}M) = \sum_{i=0}^{d} (-1)^i \mathrm{e}_Q^i(M) \binom{n+d-i}{d-i}$$
*for all* $n \geq 0$,

(iii) $QM \cap \mathrm{H}_{\mathfrak{m}}^0(M) = (0)$, *and* $Q\mathrm{H}_{\mathfrak{m}}^i(M) = (0)$ *for all* $1 \leq i \leq d-3$.

We now briefly explain how this paper is organized. In Section 2 we will summarize, for the later use in this paper, some auxiliary results on the homological degrees and torsions. We shall prove Theorem 1.3 in Section 3 (Theorem 3.4). In Section 4 we will explore examples of parameter ideals which satisfy the equality in Theorem 1.3 (1).

In what follows, unless otherwise specified, let $A$ be a Noetherian local ring with maximal ideal $\mathfrak{m}$ and $d = \dim A > 0$. Let $M$ be a finitely generated $A$-module with $s = \dim_A M$. We throughout assume that $A$ is $\mathfrak{m}$–adically complete and the field $A/\mathfrak{m}$ is infinite. For each $\mathfrak{m}$-primary ideal $I$ in $A$ we set

$$R = \mathrm{R}(I) = A[It], \quad R' = \mathrm{R}'(I) = A[It, t^{-1}], \quad \text{and} \quad \mathrm{gr}_I(A) = \mathrm{R}'(I)/t^{-1}\mathrm{R}'(I),$$

where $t$ is an indeterminate over $A$.

## 2. Preliminaries

In this section we summarize some basic properties of homological degrees and torsions of modules, which we need throughout this paper. See [GO] for the detailed proofs.

For each $j \in \mathbb{Z}$ we set

$$M_j = \mathrm{Hom}_A(\mathrm{H}_{\mathfrak{m}}^j(M), E),$$

where $E = \mathrm{E}_A(A/\mathfrak{m})$ denotes the injective envelope of $A/\mathfrak{m}$ and $\mathrm{H}_{\mathfrak{m}}^j(M)$ the $j$th local cohomology module of $M$ with respect to $\mathfrak{m}$.

We begin with the following.

**Fact 2.1.** For each $j \in \mathbb{Z}$, $M_j$ is a finitely generated $A$-module with $\dim_A M_j \leq j$, where $\dim_A(0) = -\infty$.

*Proof.* See [GO, Fact 2.1]. □

We recall the definition of homological degrees.

**Definition 2.2.** ([V2]) For each finitely generated $A$-module $M$ with $s = \dim_A M$ and for each $\mathfrak{m}$-primary ideal $I$ of $A$, we set

$$\mathrm{hdeg}_I(M) = \begin{cases} \ell_A(M) & \text{if } s \leq 0, \\ \\ \mathrm{e}_I^0(M) + \sum_{j=0}^{s-1} \binom{s-1}{j} \mathrm{hdeg}_I(M_j) & \text{if } s > 0 \end{cases}$$

and call it the homological degree of $M$ with respect to $I$.

Let us summarize some basic properties of $\mathrm{hdeg}_I(M)$.

**Fact 2.3.** Let $M$ and $M'$ are finitely generated $A$-modules. Let $I$ be an $\mathfrak{m}$-primary ideal in $A$. Then $0 \leq \mathrm{hdeg}_I(M) \in \mathbb{Z}$. We furthermore have the following:

(1) $\mathrm{hdeg}_I(M) = 0$ if and only if $M = (0)$.
(2) If $M \cong M'$, then $\mathrm{hdeg}_I(M) = \mathrm{hdeg}_I(M')$.
(3) $\mathrm{hdeg}_I(M)$ depends only on the integral closure of $I$.
(4) If $M$ is a generalized Cohen-Macaulay $A$-module, then

$$\mathrm{hdeg}_I(M) - \mathrm{e}_I^0(M) = \mathbb{I}(M)$$

and

$$\ell_A(M/QM) - \mathrm{e}_Q^0(M) \leq \mathbb{I}(M)$$

for all parameter ideals $Q$ for $M$ ([STC]), where $\mathbb{I}(M) = \sum_{j=0}^{s-1} \binom{s-1}{j} \ell_A(\mathrm{H}_\mathfrak{m}^j(M))$ denotes the Stückrad-Vogel invariant of $M$.

The following result plays a key role in the analysis of homological degree.

**Lemma 2.4.** ([V2, Proposition 3.18]) *Let* $0 \to X \to Y \to Z \to 0$ *be an exact sequence of finitely generated* $A$*-modules. Then the following assertions hold true:*

(1) *If* $\ell_A(Z) < \infty$*, then* $\mathrm{hdeg}_I(Y) \leq \mathrm{hdeg}_I(X) + \mathrm{hdeg}_I(Z)$.
(2) *If* $\ell_A(X) < \infty$*, then* $\mathrm{hdeg}_I(Y) = \mathrm{hdeg}_I(X) + \mathrm{hdeg}_I(Z)$.

*Proof.* See [GO, Lemma 2.4]. $\square$

Let $R = \mathrm{R}(I) = A[It] \subseteq A[t]$ be the Rees algebra of $I$ (here $t$ denotes an indeterminate over $A$) and let $f : I \to R,\ a \mapsto at$ be the identification of $I$ with $R_1 = It$. Set

$$\mathrm{Proj}\, R = \{\mathfrak{p} \mid \mathfrak{p} \text{ is a graded prime ideal of } R \text{ such that } \mathfrak{p} \not\supseteq R_+\}.$$

We then have the following.

**Lemma 2.5.** ([V1, Theorem 2.13]) *Let* $M$ *be a finitely generated* $A$*-module. Then there exists a finite subset* $\mathcal{F} \subseteq \mathrm{Proj}\, R$ *such that*

(1) *every* $a \in I \setminus \bigcup_{\mathfrak{p} \in \mathcal{F}} [f^{-1}(\mathfrak{p}) + \mathfrak{m}I]$ *is superficial for* $M$ *with respect to* $I$ *and*
(2) $\mathrm{hdeg}_I(M/aM) \leq \mathrm{hdeg}_I(M)$ *for each* $a \in I \setminus \bigcup_{\mathfrak{p} \in \mathcal{F}} [f^{-1}(\mathfrak{p}) + \mathfrak{m}I]$.

*Proof.* See [GO, Lemma 2.6]. $\square$

**Definition 2.6.** Let $M$ be a finitely generated $A$–module with $s = \dim_A M \geq 2$. We set

$$\mathrm{T}_I^i(M) = \sum_{j=1}^{s-i} \binom{s-i-1}{j-1} \mathrm{hdeg}_I(M_j)$$

for each $1 \leq i \leq s-1$ and call them the homological torsions of $M$ with respect to $I$.

We notice that

$$\mathrm{hdeg}_I(M) - \mathrm{T}_I^1(M) = \mathrm{e}_I^0(M) + \sum_{j=0}^{s-2} \binom{d-2}{j} \mathrm{hdeg}_I(M_i)$$

holds true. We then have the following.

**Lemma 2.7.** *Let $M$ be a finitely generated $A$–module with $s = \dim_A M \geq 3$ and $I$ an $\mathfrak{m}$-primary ideal of $A$. Then, there exists a finite subset $\mathcal{F} \subseteq \operatorname{Proj} R$ such that every $a \in I \setminus \bigcup_{\mathfrak{p}\in\mathcal{F}}[f^{-1}(\mathfrak{p}) + \mathfrak{m}I]$ is superficial for $M$ with respect to $I$, satisfying the inequality*

$$\operatorname{hdeg}_I(M/aM) - \mathrm{T}_I^1(M/aM) \leq \operatorname{hdeg}_I(M) - \mathrm{T}_I^1(M).$$

*Proof.* Thanks to Lemma 2.5, there exists a finite subset $\mathcal{F} \subseteq \operatorname{Proj} R$ such that every $a \in I \setminus \bigcup_{\mathfrak{p}\in\mathcal{F}}[f^{-1}(\mathfrak{p}) + \mathfrak{m}I]$ is superficial for $M$ and $M_j$ with respect to $I$ and $\operatorname{hdeg}_I(M_j/aM_j) \leq \operatorname{hdeg}_I(M_j)$ for all $1 \leq j \leq s-2$. Set $\overline{M} = M/aM$. Consider the long exact sequence

$$0 \to (0) :_M a \to \mathrm{H}^0_{\mathfrak{m}}(M) \xrightarrow{a} \mathrm{H}^0_{\mathfrak{m}}(M) \to \mathrm{H}^0_{\mathfrak{m}}(\overline{M}) \to \mathrm{H}^1_{\mathfrak{m}}(M) \xrightarrow{a} \mathrm{H}^1_{\mathfrak{m}}(M) \to \mathrm{H}^1_{\mathfrak{m}}(\overline{M}) \to \cdots$$

$$\cdots \to \mathrm{H}^j_{\mathfrak{m}}(M) \xrightarrow{a} \mathrm{H}^j_{\mathfrak{m}}(M) \to \mathrm{H}^j_{\mathfrak{m}}(\overline{M}) \to \mathrm{H}^{j+1}_{\mathfrak{m}}(M) \xrightarrow{a} \mathrm{H}^{j+1}_{\mathfrak{m}}(M) \to \cdots$$

of local cohomology modules induced from the exact sequence

$$0 \to (0) :_M a \to M \xrightarrow{a} M \to \overline{M} \to 0.$$

Then, taking the Matlis dual of the above long exact sequence, we get exact sequences

$$0 \to M_{j+1}/aM_{j+1} \to \overline{M}_j \to (0) :_{M_j} a \to 0$$

and embeddings

$$0 \to (0) :_{M_j} a \to M_j$$

for all $0 \leq j \leq s-3$. Consequently, because $\ell_A((0) :_{M_j} a) < \infty$ and by Lemma 2.4 we have

$$\begin{aligned} \operatorname{hdeg}_I(\overline{M}_j) &\leq \operatorname{hdeg}_I([(0) :_{M_j} a]) + \operatorname{hdeg}_I(M_{j+1}/aM_{j+1}) \\ &\leq \operatorname{hdeg}_I(M_j) + \operatorname{hdeg}_I(M_{j+1}) \end{aligned}$$

for each $0 \leq j \leq s-3$. Hence, because $\mathrm{e}^0_Q(M) = \mathrm{e}^0_Q(\overline{M})$, we get

$$\begin{aligned} \operatorname{hdeg}_I(\overline{M}) - \mathrm{T}_I^1(\overline{M}) &= \mathrm{e}^0_I(\overline{M}) + \sum_{j=0}^{s-3} \binom{s-3}{j} \operatorname{hdeg}_I(\overline{M}_j) \\ &\leq \mathrm{e}^0_I(\overline{M}) + \sum_{j=0}^{s-3} \binom{s-3}{j} \{\operatorname{hdeg}_I(M_j) + \operatorname{hdeg}_I(M_{j+1})\} \\ &= \mathrm{e}^0_I(M) + \sum_{j=0}^{s-2} \binom{s-2}{j} \operatorname{hdeg}_I(M_j) \\ &= \operatorname{hdeg}_I(M) - \mathrm{T}_I^1(M), \end{aligned}$$

as required. $\square$

## 3. The sectional genera and the homological degrees of parameters

In this section we study the behavior of the sectional genera of parameters. Let $Q = (a_1, a_2, \ldots, a_s)$ be a parameter ideal for $M$. We set

$$\mathrm{g}_s(Q; M) = \ell_A(M/QM) - \mathrm{e}_Q^0(M) + \mathrm{e}_Q^1(M)$$

and call it the sectional genus of $M$ with respect to $Q$.

We denote by $\mathrm{H}_i(Q; M)$ ($i \in \mathbb{Z}$) the $i$–th homology module of the Koszul complex $\mathrm{K}_\bullet(Q; M)$ generated by the system $a_1, a_2, \ldots, a_s$ of parameters for $M$. Set

$$\chi_1(Q; M) = \sum_{i \geq 1} (-1)^{i-1} \ell_A(\mathrm{H}_i(Q; M))$$

and call it the first Euler characteristic of $M$ relative to $Q$. Hence

$$\chi_1(Q; M) = \ell_A(M/QM) - \mathrm{e}_Q^0(M) \geq 0$$

and so that

$$\mathrm{g}_s(Q; M) = \chi_1(Q; M) + \mathrm{e}_Q^1(M).$$

The following result is due to [GHV]. We indicate a brief proof for the sake of completeness.

**Lemma 3.1.** ([GHV, Proposition 4.1]) *Let $M$ be a finitely generated $A$-module with $\dim_A M = 1$. Let $Q = (a)$ be a parameter ideal for $M$. Then*

$$\mathrm{g}_s(Q; M) \leq 0$$

*and the following two conditions are equivalent:*

(1) $\mathrm{g}_s(Q; M) = 0$,

(2) *$a$ forms a $d$-sequence on $M$.*

*Proof.* We notice that $\mathrm{e}_Q^1(M) = -\ell_A(\mathrm{H}_\mathfrak{m}^0(M))$ (see [GNi, Lemma 2.4], [MSV, Proposition 3.1]) and $\chi_1(Q; M) = \ell_A([(0) :_M a]) \leq \ell_A(\mathrm{H}_\mathfrak{m}^0(M))$. Therefore we get $\mathrm{g}_s(Q; M) = \chi_1(Q; M) + \mathrm{e}_Q^1(M) \leq 0$ as required.

Let us consider the second assertion.

(1) $\Rightarrow$ (2) We have $\ell_A([(0) :_M a]) = \chi_1(Q; M) = -\mathrm{e}_Q^1(M) = \ell_A(\mathrm{H}_\mathfrak{m}^0(M))$ so that $[(0) :_M a^2] \subseteq \mathrm{H}_\mathfrak{m}^0(M) = [(0) :_M a]$. Thus $a$ forms a $d$-sequence on $M$.

(2) $\Rightarrow$ (1) Since $a$ forms a $d$-sequence on $M$, we have $\mathrm{H}_\mathfrak{m}^0(M) = [(0) :_M a]$. Hence we get $\mathrm{g}_s(Q; M) = \chi_1(Q; M) + \mathrm{e}_Q^1(M) = 0$ as required. $\square$

We note the following.

**Lemma 3.2.** *Let $M$ be a finitely generated $A$-module with $s = \dim_A M \geq 2$. Let $Q$ be a parameter ideal for $M$ and assume that $a \in Q \backslash \mathfrak{m}Q$ is a superficial element for $M$ with respect to $Q$. Then $\mathrm{g}_s(Q; M) = \mathrm{g}_s(\overline{Q}; \overline{M}) + \ell_A([(0) :_M a])$ if $d = 2$, and $\mathrm{g}_s(Q; M) = \mathrm{g}_s(\overline{Q}; \overline{M})$ if $d \geq 3$, where $\overline{M} = M/aM$ and $\overline{Q} = Q/(a)$.*

*Proof.* We have $\ell_A(M/QM) = \ell_A(\overline{M}/Q\overline{M})$ and $\mathrm{e}_Q^0(M) = \mathrm{e}_Q^0(\overline{M})$. We also have $\mathrm{e}_Q^1(M) = \mathrm{e}_Q^1(\overline{M}) + \ell_A([(0) :_M a])$ if $d = 2$, and $\mathrm{e}_Q^1(M) = \mathrm{e}_Q^1(\overline{M})$ if $d \geq 3$ ([N, (22.6)]). Thus we get the rquired equalities. $\square$

The following result gives an upper bound for $\mathrm{g}_s(Q;M)$.

**Proposition 3.3.** *Let $M$ be a finitely generated $A$-module with $d = \dim_A M \geq 2$. Then*

$$\mathrm{g}_s(Q;M) \leq \mathrm{hdeg}_Q(M) - \mathrm{e}_Q^0(M) - \mathrm{T}_Q^1(M)$$

*for every parameter ideal $Q$ of $A$.*

*Proof.* Suppose $d = 2$. We choose an element $a \in Q \backslash \mathfrak{m}Q$ so that $a$ is superficial for $M$ with respect to $Q$. Let $\overline{M} = M/aM$ and $\overline{Q} = Q/(a)$. Then we have

$$\mathrm{g}_s(Q;M) = \mathrm{g}_s(\overline{Q};\overline{M}) + \ell_A([(0):_M a]) \leq \ell_A([(0):_M a]) \leq \ell_A(\mathrm{H}_{\mathfrak{m}}^0(M))$$

by Lemma 3.1 and Lemma 3.2.

Assume that $d \geq 3$ and that our assertion holds true for $d-1$. We choose an element $a \in Q \backslash \mathfrak{m}Q$ so that $a$ is superficial for $M$ with respect to $Q$ and $\mathrm{hdeg}_Q(M/aM) - \mathrm{T}_Q^1(M/aM) \leq \mathrm{hdeg}_Q(M) - \mathrm{T}_Q^1(M)$ (Lemma 2.7). Then, setting $\overline{M} = M/aM$ and $\overline{Q} = Q/(a)$, by the hypothesis of induction on $d$ we get

$$\begin{aligned} \mathrm{g}_s(Q;M) = \mathrm{g}_s(\overline{Q};\overline{M}) &\leq \mathrm{hdeg}_Q(\overline{M}) - \mathrm{e}_Q^0(\overline{M}) - \mathrm{T}_Q^1(\overline{M}) \\ &\leq \mathrm{hdeg}_Q(M) - \mathrm{e}_Q^0(M) - \mathrm{T}_Q^1(M) \end{aligned}$$

by Lemma 3.2. □

We notice here that, in [GhGHOPV], it was proved that, for parameter ideals $Q$ for $M$, the upper bound

$$\chi_1(Q;M) \leq \mathrm{hdeg}_Q(M) - \mathrm{e}_Q^0(M)$$

of $\chi_1(Q;M)$ and the lower bound

$$\mathrm{e}_Q^1(M) \geq -\mathrm{T}_Q^1(M)$$

of $\mathrm{e}_Q^1(M)$. Thus the upper bound for $\mathrm{g}_s(Q;M) = \chi_1(Q;M) + \mathrm{e}_Q^1(M)$ is given by the sum of the upper bound $\mathrm{hdeg}_Q(M) - \mathrm{e}_Q^0(M)$ of $\chi_1(Q;M)$ and the lower bound $-\mathrm{T}_Q^1(M)$ of $\mathrm{e}_Q^1(M)$. It seems natural to ask what happens on the parameters $Q$ for $M$, when $\mathrm{g}_s(Q;M) = \mathrm{hdeg}_Q(M) - \mathrm{e}_Q^0(M) - \mathrm{T}_Q^1(M)$. The following theorem answers the question, which is the main result of this paper (Theorem 1.3).

**Theorem 3.4.** *Let $M$ be a finitely generated $A$-module with $d = \dim_A M \geq 2$. Let $Q$ be a parameter ideal of $A$. Then the following conditions are equivalent:*

(1) $\mathrm{g}_s(Q;M) = \mathrm{hdeg}_Q(M) - \mathrm{e}_Q^0(M) - \mathrm{T}_Q^1(M)$.

(2) *The following conditions are satisfied:*

(a)

$$(-1)^i \mathrm{e}_Q^i(M) = \begin{cases} \mathrm{T}_Q^i(M) & \text{if } 2 \leq i \leq d-1, \\ \ell_A(\mathrm{H}_{\mathfrak{m}}^0(M)) & \text{if } i = d \end{cases}$$

*for all $2 \leq i \leq d$ and*

(b) $\ell_A(M/QM) - \sum_{i=0}^{d} (-1)^i \mathrm{e}_Q^i(M) = 0$.

*When this is the case, we have the following:*

(i) *there exist elements* $a_1, a_2, \ldots, a_d \in A$ *such that* $Q = (a_1, a_2, \ldots, a_d)$ *and* $a_1, a_2, \ldots, a_d$ *forms a* $d$*-sequence on* $M$,

(ii)
$$\ell_A(M/Q^{n+1}M) = \sum_{i=0}^{d}(-1)^i \mathrm{e}_Q^i(M)\binom{n+d-i}{d-i}$$
*for all* $n \geq 0$,

(iii) $QM \cap \mathrm{H}_{\mathfrak{m}}^0(M) = (0)$ *and* $Q\mathrm{H}_{\mathfrak{m}}^i(M) = (0)$ *for all* $1 \leq i \leq d-3$.

To prove Theorem 3.4, we need the following:

**Lemma 3.5.** *Let* $M$ *be a finitely generated* $A$*-module with* $d = \dim_A M \geq 2$ *and let* $Q$ *be a parameter ideal of* $A$*. Then* $\mathrm{g}_s(Q; M) = \mathrm{hdeg}_Q(M) - \mathrm{e}_Q^0(M) - \mathrm{T}_Q^1(M)$ *if and only if* $\mathrm{g}_s(Q; M/\mathrm{H}_{\mathfrak{m}}^0(M)) = \mathrm{hdeg}_Q(M/\mathrm{H}_{\mathfrak{m}}^0(M)) - \mathrm{e}_Q^0(M/\mathrm{H}_{\mathfrak{m}}^0(M)) - \mathrm{T}_Q^1(M/\mathrm{H}_{\mathfrak{m}}^0(M))$ *and* $QM \cap \mathrm{H}_{\mathfrak{m}}^0(M) = (0)$.

*Proof.* We set $W = \mathrm{H}_{\mathfrak{m}}^0(M)$ and $M' = M/W$. Consider the exact sequence
$$0 \to W/[QM \cap W] \to M/QM \to M'/QM' \to 0 \qquad (\sharp)$$
obtained by the canonical exact sequence
$$0 \to W \to M \to M' \to 0.$$
Assume that $\mathrm{g}_s(Q; M) = \mathrm{hdeg}_Q(M) - \mathrm{e}_Q^0(M) - \mathrm{T}_Q^1(M)$. Then we have
$$\begin{aligned}
\mathrm{g}_s(Q; M) &= \ell_A(M/QM) - \mathrm{e}_Q^0(M) + \mathrm{e}_Q^1(M)\\
&= \{\ell_A(M'/QM') + \ell_A(W/[QM \cap W])\} - \mathrm{e}_Q^0(M') + \mathrm{e}_Q^1(M')\\
&= \mathrm{g}_s(Q; M') + \{\ell_A(W) - \ell_A(QM \cap W)\}\\
&\leq \mathrm{hdeg}_Q(M') - \mathrm{e}_Q^0(M') - \mathrm{T}_Q^1(M') + \ell_A(W)\\
&= \mathrm{hdeg}_Q(M) - \mathrm{e}_Q^0(M) - \mathrm{T}_Q^1(M) = \mathrm{g}_s(Q; M),
\end{aligned}$$
because $\mathrm{e}_Q^0(M) = \mathrm{e}_Q^0(M')$, $\mathrm{e}_Q^1(M) = \mathrm{e}_Q^1(M')$, $\ell_A(M/QM) = \ell_A(M'/QM') + \ell_A(W/[QM \cap W])$ by exact sequence $(\sharp)$, $\mathrm{g}_s(Q; M') \leq \mathrm{hdeg}_Q(M') - \mathrm{e}_Q^0(M') - \mathrm{T}_Q^1(M')$ by Proposition 3.3, $\mathrm{hdeg}_Q(M) = \mathrm{hdeg}_Q(M') + \ell_A(W)$ by Lemma 2.4, and $\mathrm{T}_Q^1(M) = \mathrm{T}_Q^1(M')$. Thus $\mathrm{g}_s(Q; M') = \mathrm{hdeg}_Q(M') - \mathrm{e}_Q^0(M') - \mathrm{T}_Q^1(M')$ and $QM \cap W = (0)$.

Conversely, assume that $\mathrm{g}_s(Q; M') = \mathrm{hdeg}_Q(M') - \mathrm{e}_Q^0(M') - \mathrm{T}_Q^1(M')$ and $QM \cap W = (0)$. Then
$$\begin{aligned}
\mathrm{g}_s(Q; M) &= \ell_A(M/QM) - \mathrm{e}_Q^0(M) + \mathrm{e}_Q^1(M)\\
&= \{\ell_A(M'/QM') + \ell_A(W)\} - \mathrm{e}_Q^0(M') + \mathrm{e}_Q^1(M')\\
&= \mathrm{g}_s(Q; M') + \ell_A(W)\\
&= \mathrm{hdeg}_Q(M') - \mathrm{e}_Q^0(M') - \mathrm{T}_Q^1(M') + \ell_A(W)\\
&= \mathrm{hdeg}_Q(M) - \mathrm{e}_Q^0(M) - \mathrm{T}_Q^1(M),
\end{aligned}$$

because $e^0_Q(M) = e^0_Q(M')$, $e^1_Q(M) = e^1_Q(M')$, $\ell_A(M/QM) = \ell_A(M'/QM') + \ell_A(W)$ by exact sequence $(\sharp)$, $\mathrm{hdeg}_Q(M) = \mathrm{hdeg}_Q(M') + \ell_A(W)$ by Lemma 2.4, and $\mathrm{T}^1_Q(M) = \mathrm{T}^1_Q(M')$ which proves Lemma 3.5. $\square$

The following result shows that Theorem 3.4 (i) (Theorem 1.3 (i)) holds true, once $\mathrm{g}_s(Q;M) = \mathrm{hdeg}_Q(M) - e^0_Q(M) - \mathrm{T}^1_Q(M)$.

**Proposition 3.6.** *Let $M$ be a finitely generated $A$-module with $d = \dim_A M \geq 2$ and $Q$ a parameter ideal of $A$. Let $a_1 \in Q \backslash \mathfrak{m}Q$ be a superficial element for $M$ with respect to $Q$ such that $\mathrm{hdeg}_Q(M/a_1M) - \mathrm{T}^1_Q(M/a_1M) \leq \mathrm{hdeg}_Q(M) - \mathrm{T}^1_Q(M)$. Assume that*

$$\mathrm{g}_s(Q;M) = \mathrm{hdeg}_Q(M) - e^0_Q(M) - \mathrm{T}^1_Q(M).$$

*Then there exist elements $a_2, a_3, \ldots, a_d \in A$ such that $Q = (a_1, a_2, \ldots, a_d)$ and $a_1, a_2, \ldots, a_d$ forms a $d$-sequence on $M$.*

We note the following Lemma 3.7, before giving a proof of Proposition 3.6. The following result is, more or less, known. Let us indicate a brief proof for the sake of completeness, because it plays a key role in our proof of Proposition 3.6.

**Lemma 3.7.** *Let $M$ be a finitely generated $A$-module and $n > 0$ an integer. Let $a_1, a_2, \ldots, a_n \in A$ and assume that $a_1$ is a superficial element for $M$. Then $a_1, a_2, \ldots, a_n$ forms a $d$-sequence on $M$ if and only if $a_2, a_3, \ldots, a_n$ forms a $d$-sequence on $M/a_1M$ and $(a_1, a_2, \ldots, a_n)M \cap \mathrm{H}^0_{\mathfrak{m}}(M) = (0)$.*

*Proof.* Assume that $a_2, a_3, \ldots, a_n$ forms a $d$-sequence on $M/a_1M$ and $(a_1, a_2, \ldots, a_n)M \cap \mathrm{H}^0_{\mathfrak{m}}(M) = (0)$, and set $\overline{M} = M/a_1M$. Then we have

$$[(a_2, a_3, \ldots, a_{i-1})\overline{M} :_{\overline{M}} a_i a_j] = [(a_2, a_3, \ldots, a_{i-1})\overline{M} :_{\overline{M}} a_j],$$

so that

$$[(a_1, a_2, \ldots, a_{i-1})M :_M a_i a_j] = [(a_1, a_2, \ldots, a_{i-1})M :_M a_j]$$

for all $2 \leq i \leq j \leq n$. It is now enough to show that $[(0) :_M a_1 a_j] = [(0) :_M a_j]$ for all $1 \leq j \leq n$. Take $m \in [(0) :_M a_1 a_j]$. Then $a_1 a_j m = 0$. Then

$$a_j m \in [(0) :_M a_1] \cap (a_1, a_2, \ldots, a_n)M \subseteq W \cap (a_1, a_2, \ldots, a_n)M = (0),$$

because $a_1$ is superficial for $M$. Hence $m \in [(0) :_M a_j]$, so that $[(0) :_M a_1 a_j] \subseteq [(0) :_M a_j]$. Thus $[(0) :_M a_1 a_j] = [(0) :_M a_j]$ for all $1 \leq j \leq n$. Hence $a_1, a_2, \ldots, a_n$ forms a $d$-sequence on $M$. The converse holds true by the definition of a $d$-sequence. This completes the proof of Proposition 3.7. $\square$

*Proof of Proposition 3.6.* We proceed by induction on $d$. Set $\overline{M} = M/a_1M$, $\overline{A} = A/(a_1)$, and $\overline{Q} = Q/(a_1)$. Suppose that $d = 2$. Let $Q = (a_1, a_2)$. Then we have

$$\mathrm{g}_s(Q;M) \;\; = \;\; \mathrm{g}_s(\overline{Q};\overline{M}) + \ell_A([(0) :_M a_1]) \leq \mathrm{hdeg}_Q(M_0) = \mathrm{g}_s(Q;M)$$

because $\mathrm{g}_s(Q;M) = \mathrm{g}_s(\overline{Q};\overline{M}) + \ell_A([(0) :_M a_1])$ by Lemma 3.2 and $\mathrm{g}_s(\overline{Q};\overline{M}) \leq 0$ by Lemma 3.1. Hence $\mathrm{g}_s(\overline{Q};\overline{M}) = 0$ so that $a_2$ forms a $d$-sequence on $\overline{M}$ by Lemma 3.1. Therefore, since $QM \cap \mathrm{H}^0_{\mathfrak{m}}(M) = (0)$ by Lemma 3.5, $a_1, a_2$ forms a $d$-sequence on $M$ by Lemma 3.7.

Assume that $d \geq 3$ and that our assertion holds true for $d-1$. Then since

$$\begin{aligned} \mathrm{g}_s(Q;M) = \mathrm{g}_s(\overline{Q};\overline{M}) &\leq \mathrm{hdeg}_Q(\overline{M}) - \mathrm{e}_Q^0(\overline{M}) - \mathrm{T}_Q^1(\overline{M}) \\ &\leq \mathrm{hdeg}_Q(M) - \mathrm{e}_Q^0(M) - \mathrm{T}_Q^1(M) \\ &= \mathrm{g}_s(Q;M) \end{aligned}$$

by Lemma 3.2 and Proposition 3.3, we have $\mathrm{g}_s(\overline{Q};\overline{M}) = \mathrm{hdeg}_Q(\overline{M}) - \mathrm{e}_Q^0(\overline{M}) - \mathrm{T}_Q^1(\overline{M})$. Because the residue class field $A/\mathfrak{m}$ of $A$ is infinite, we may also choose an element $a_2 \in Q$ so that $a_2$ is superficial for $\overline{M}$ with respect to $\overline{Q}$, $\mathrm{hdeg}_Q(\overline{M}/a_2\overline{M}) - \mathrm{T}_Q^1(\overline{M}/a_2\overline{M}) \leq \mathrm{hdeg}_Q(\overline{M}) - \mathrm{T}_Q^1(\overline{M})$ (Lemma 2.7), and $a_1, a_2$ forms, furthermore, a part of a minimal system of generators of $Q$. Then the hypothesis of induction on $d$ guarantees that there exist elements $a_3, a_4, \ldots, a_d \in A$ such that $\overline{Q} = (a_2, a_3, \ldots, a_d)\overline{A}$ and $a_2, a_3, \ldots, a_d$ forms a $d$-sequence on $\overline{M}$. Thus, thanks to Lemma 3.7, $a_1, a_2, \cdots, a_d$ forms a $d$-sequence on $M$ because $QM \cap \mathrm{H}_\mathfrak{m}^0(M) = (0)$ by Lemma 3.5. □

The following result plays a key role in our proof of our main theorem.

**Proposition 3.8.** (*cf.* [GO, Proposition 3.7]) *Let $M$ be a finitely generated $A$-module with $d = \dim_A M$. Let $Q = (a_1, a_2, \ldots, a_d)$ be a parameter ideal of $A$ and assume that $a_1, a_2, \ldots, a_d$ forms a $d$-sequence on $M$. Then we have the following, where $Q_i = (a_1, a_2, \ldots, a_i)$ for each $0 \leq i \leq d$.*

1. $\mathrm{e}_Q^0(M) = \ell_A(M/QM) - \ell_A([Q_{d-1}M :_M a_d]/Q_{d-1}M)$.
2. $(-1)^i\mathrm{e}_Q^i(M) = \ell_A(\mathrm{H}_\mathfrak{m}^0(M/Q_{d-i}M)) - \ell_A(\mathrm{H}_\mathfrak{m}^0(M/Q_{d-i-1}M))$ *for* $1 \leq i \leq d-1$ *and* $(-1)^d\mathrm{e}_Q^d(M) = \ell_A(\mathrm{H}_\mathfrak{m}^0(M))$.
3. $\ell_A(M/Q^{n+1}M) = \sum_{i=0}^d (-1)^i \mathrm{e}_Q^i(M)\binom{n+d-i}{d-i}$ *for all* $n \geq 0$. *Hence* $\ell_A(M/QM) - \sum_{i=0}^d (-1)^i\mathrm{e}_Q^i(M) = 0$.

We are now in a position to prove Theorem 3.4 (Theorem 1.3).

*Proof of Theorem 3.4.* (1) ⇒ (2) Since the last assertions (i) and (ii) follow from Proposition 3.6 and Proposition 3.8, we have assertion (b). It is now enough to show that assertion (a) holds true. We proceed by induction on $d$. Suppose that $d = 2$. Then, because $a_1, a_2$ forms a $d$-sequence on $M$ by Proposition 3.6, we have $\mathrm{e}_Q^2(M) = \ell_A(\mathrm{H}_\mathfrak{m}^0(M))$ by Proposition 3.8.

Assume that $d \geq 3$ and that our assertion holds true for $d-1$. Choose an element $a \in Q \backslash \mathfrak{m}Q$ so that $a$ is superficial for $M$ and $M_j$ with respect to $Q$ and $\mathrm{hdeg}_Q(M_j/aM_j) \leq \mathrm{hdeg}_Q(M_j)$ for all $1 \leq j \leq d-2$ (Lemma 2.5), and set $\overline{M} = M/aM$ and $\overline{Q} = Q/(a)$. Then by the same argument as is in the proof of Lemma 2.7, we get the inequalities

$$\mathrm{hdeg}_Q(\overline{M}_j) \leq \mathrm{hdeg}_Q([(0) :_{M_j} a]) + \mathrm{hdeg}_Q(M_{j+1}/aM_{j+1})$$

and

$$\ell_A([(0) :_{M_j} a]) \leq \mathrm{hdeg}_Q(M_j)$$

for all $0 \le j \le d-3$. Hence

$$\begin{aligned}
\mathrm{g}_s(Q;M) = \mathrm{g}_s(\overline{Q};\overline{M}) &\le \mathrm{hdeg}_Q(\overline{M}) - \mathrm{e}_Q^0(\overline{M}) - \mathrm{T}_Q^1(\overline{M}) \\
&= \sum_{j=0}^{d-3} \binom{d-3}{j} \mathrm{hdeg}_Q(\overline{M}_j) \\
&\le \sum_{j=0}^{d-3} \binom{d-3}{j} \{\mathrm{hdeg}_Q([(0):_{M_j} a]) + \mathrm{hdeg}_Q(M_{j+1}/aM_{j+1})\} \\
&\le \sum_{j=0}^{d-3} \binom{d-3}{j} \{\mathrm{hdeg}_Q(M_j) + \mathrm{hdeg}_Q(M_{j+1})\} \\
&= \sum_{j=0}^{d-2} \binom{d-2}{j} \mathrm{hdeg}_Q(M_j) \\
&= \mathrm{hdeg}_Q(M) - \mathrm{e}_Q^0(M) - \mathrm{T}_Q^1(M) = \mathrm{g}_s(Q;M),
\end{aligned}$$

because $\mathrm{g}_s(Q;M) = \mathrm{g}_s(\overline{Q};\overline{M})$ by Lemma 3.2 and $\mathrm{g}_s(\overline{Q};\overline{M}) \le \mathrm{hdeg}_Q(\overline{M}) - \mathrm{e}_Q^0(\overline{M}) - \mathrm{T}_Q^1(\overline{M})$ by Proposition 3.3. Thus

$$\mathrm{g}_s(\overline{Q};\overline{M}) = \mathrm{hdeg}_Q(\overline{M}) - \mathrm{e}_Q^0(\overline{M}) - \mathrm{T}_Q^1(\overline{M}),$$

$$\mathrm{hdeg}_Q(\overline{M}_j) = \mathrm{hdeg}_Q(M_j) + \mathrm{hdeg}_Q(M_{j+1}),$$

and $aM_j = (0)$ for all $0 \le j \le d-3$. On the other hand, since $a$ is superficial for $M$ with respect to $Q$, we have $\mathrm{e}_Q^i(M) = \mathrm{e}_Q^i(\overline{M})$ for all $0 \le i \le d-2$ and $(-1)^{d-1}\mathrm{e}_Q^{d-1}(M) = (-1)^{d-1}\mathrm{e}_Q^{d-1}(\overline{M}) - \ell_A([(0):_M a])$ ([N, (22.6)]). Therefore the hypothesis of induction on $d$ yields that

$$\begin{aligned}
(-1)^i\mathrm{e}_Q^i(M) = (-1)^i\mathrm{e}_Q^i(\overline{M}) &= \mathrm{T}_Q^i(\overline{M}) \\
&= \sum_{j=1}^{d-1-i} \binom{d-i-2}{j-1} \mathrm{hdeg}_Q(\overline{M}_j) \\
&= \sum_{j=1}^{d-1-i} \binom{d-i-2}{j-1} \{\mathrm{hdeg}_Q(M_j) + \mathrm{hdeg}_Q(M_{j+1})\} \\
&= \sum_{j=1}^{d-i} \binom{d-i-1}{j-1} \mathrm{hdeg}_Q(M_j) \\
&= \mathrm{T}_Q^i(M)
\end{aligned}$$

for $1 \le i \le d-2$ and that

$$\begin{aligned}
(-1)^{d-1}\mathrm{e}_Q^{d-1}(M) &= (-1)^{d-1}\mathrm{e}_Q^{d-1}(\overline{M}) - \ell_A([(0):_M a]) \\
&= \ell_A(\mathrm{H}_{\mathfrak{m}}^0(\overline{M})) - \ell_A(\mathrm{H}_{\mathfrak{m}}^0(M)) \\
&= \{\mathrm{hdeg}_Q(M_0) + \mathrm{hdeg}_Q(M_1)\} - \mathrm{hdeg}_Q(M_0) \\
&= \mathrm{hdeg}_Q(M_1) \\
&= \mathrm{T}_Q^{d-1}(M),
\end{aligned}$$

because $a\mathrm{H}_{\mathfrak{m}}^0(M) = (0)$ and $\ell_A(\mathrm{H}_{\mathfrak{m}}^0(\overline{M})) = \mathrm{hdeg}_Q(\overline{M}_0) = \mathrm{hdeg}_Q(M_0) + \mathrm{hdeg}_Q(M_1)$. Thus, as the equality $(-1)^d\mathrm{e}_Q^d(M) = \ell_A(\mathrm{H}_{\mathfrak{m}}^0(M))$ holds true by Proposition 3.8, assertion (a) follows, which proves the implication (1) $\Rightarrow$ (2).

(2) $\Rightarrow$ (1) We have

$$\begin{aligned}
\sum_{i=2}^{d-1}\mathrm{T}_Q^i(M) &= \sum_{i=2}^{d-1}\sum_{j=1}^{d-i}\binom{d-i-1}{j-1}\mathrm{hdeg}_Q(M_j) \\
&= \sum_{j=1}^{d-2}\left\{\sum_{i=2}^{d-j}\binom{d-i-1}{j-1}\right\}\mathrm{hdeg}_Q(M_j) \\
&= \sum_{j=1}^{d-2}\left\{\sum_{i=2}^{d-j}\left[\binom{d-i}{j} - \binom{d-i-1}{j}\right]\right\}\mathrm{hdeg}_Q(M_j) \\
&= \sum_{j=1}^{d-2}\left\{\sum_{i=2}^{d-j}\binom{d-i}{j} - \sum_{i=2}^{d-j-1}\binom{d-i-1}{j}\right\}\mathrm{hdeg}_Q(M_j) \\
&= \sum_{j=1}^{d-2}\left\{\sum_{i=2}^{d-j}\binom{d-i}{j} - \sum_{i=3}^{d-j}\binom{d-i}{j}\right\}\mathrm{hdeg}_Q(M_j) \\
&= \sum_{j=1}^{d-2}\binom{d-2}{j}\mathrm{hdeg}_Q(M_j).
\end{aligned}$$

Thus

$$\begin{aligned}
\mathrm{g}_s(Q;M) &= \ell_A(M/QM) - \mathrm{e}_Q^0(M) + \mathrm{e}_Q^1(M) \\
&= \sum_{i=2}^{d}(-1)^i\mathrm{e}_Q^i(M) \\
&= \sum_{i=2}^{d-1}\mathrm{T}_Q^i(M) + \ell_A(\mathrm{H}_{\mathfrak{m}}^0(M)) \\
&= \sum_{j=0}^{d-2}\binom{d-2}{j}\mathrm{hdeg}_Q(M_j) = \mathrm{hdeg}_Q(M) - \mathrm{e}_Q^0(M) - \mathrm{T}_Q^1(M),
\end{aligned}$$

which shows the implication (2) $\Rightarrow$ (1).

We now consider assertion (iii). We get $QM \cap \mathrm{H}^0_{\mathfrak{m}}(M) = (0)$ by Lemma 3.5. Suppose that $d \geq 4$. Let $Q = (a_1, a_2, \ldots, a_d)$ and $1 \leq i \leq d$. Since the residue class field $A/\mathfrak{m}$ of $A$ is infinite, we may choose the elements $a_i's$ so that $a_i$ is superficial for $M$ and $M_j$ with respect to $Q$ and $\mathrm{hdeg}_Q(M_j/a_iM_j) \leq \mathrm{hdeg}_Q(M_j)$ for all $1 \leq j \leq d-2$. Then, thanks to the proof of the implication (1) $\Rightarrow$ (2), $a_iM_j = (0)$ for all $1 \leq j \leq d-3$. Consequently, by the symmetry of $a_i's$, $Q\mathrm{H}^j_{\mathfrak{m}}(M) = (0)$ for all $1 \leq j \leq d-3$, which proves assertion (iii) and Theorem 3.4. $\square$

## 4. Examples

In this section we will explore examples of parameter ideals $Q$ which satisfy the equality $\mathrm{g}_s(Q; A) = \mathrm{hdeg}_Q(A) - \mathrm{e}^0_Q(A) - \mathrm{T}^1_Q(A)$ but $A$ is not a generalized Cohen-Macaulay local ring.

Let us begin with the following.

**Proposition 4.1.** *Let $R$ be a Cohen-Macaulay local ring with maximal ideal $\mathfrak{n}$ and $d = \dim R > 0$. Let $M$ be a Cohen-Macaulay $R$-module with $\dim_R M = d-1$. We set $A = R \ltimes M$ be an idealization of $M$ over $R$. Let $J$ be a parameter ideal in $R$ and $Q = JA$. Then we have the following.*

1. *$\dim A = d$ and $\mathrm{depth} A = d-1$, but $\mathrm{H}^{d-1}_{\mathfrak{n}}(M)$ is not finitely generated.*
2. *$\mathrm{g}_s(Q; A) = \ell_R(M/JM) - \mathrm{e}^0_J(M)$.*
3. *$\mathrm{hdeg}_Q(A) - \mathrm{e}^0_Q(A) + \mathrm{T}^1_Q(A) = 0$.*
4. *Therefore $\mathrm{g}_s(Q; A) = \mathrm{hdeg}_Q(A) - \mathrm{e}^0_Q(A) - \mathrm{T}^1_Q(A)$ if and only if $\ell_R(M/JM) = \mathrm{e}^0_J(M)$.*

*Proof.* Since $\mathrm{depth} A = d-1$, $\mathrm{hdeg}_Q(A) - \mathrm{e}^0_Q(A) + \mathrm{T}^1_Q(A) = 0$ holds true. We have

$$\begin{aligned} \ell_A(A/Q^{n+1}) &= \ell_R(R/J^{n+1}) + \ell_R(M/J^{n+1}M) \\ &= \ell_R(R/J)\binom{n+d}{d} + \left\{\sum_{i=0}^{d-1}(-1)^i \mathrm{e}^i_J(M)\binom{n+d-1-i}{d-1-i}\right\} \end{aligned}$$

for all $n \gg 0$. Therefore $\mathrm{e}^0_Q(A) = \ell_R(R/J)$ and $\mathrm{e}^1_Q(A) = -\mathrm{e}^0_J(M)$ so that $\mathrm{g}_s(Q; A) = \ell_A(A/Q) - \mathrm{e}^0_Q(A) + \mathrm{e}^1_Q(A) = \ell_R(M/JM) - \mathrm{e}^0_J(M)$ as required. $\square$

The following example shows that there exists a parameter ideal for a Cohen-Macaulay module which satisfies the condition of Proposition 4.1 (4), where the finitely generated module $N$ over a Cohen-Macaulay local ring $R$ is said to be an Ulrich module with respect to an $\mathfrak{m}$-primary ideal $I$, if the following three conditions are satisfied ([GOTWY]).

1. $N$ is a maximal Cohen-Macaulay $R$-module, that is $\mathrm{depth}_R M = \dim R$,
2. $\mathrm{e}^0_I(N) = \ell_R(N/IN)$, and
3. $N/IN$ is $R/I$-free.

**Example 4.2.** Let $R$ be a Cohen-Macaulay local ring with $d = \dim R \geq 1$. Let $I = (x_1, x_2, \cdots, x_d)$ be a parameter ideal of $R$. Let

$$0 \to R^d \xrightarrow{\partial} R^d \to C \to 0 \qquad (\flat)$$

be the exact sequence of $R$-modules where the $d \times d$ matrix $\partial$ has the form $\partial = (\partial_{ij})_{1 \le i,j \le d}$ with

$$\partial_{ij} = \begin{cases} x_1 & \text{if } i = j, \\ x_{j-i+1} & \text{if } i < j \\ 0 & \text{if } i > j \end{cases}$$

and $C = \mathrm{Coker}\partial$. Then we have the following.

(1) $C$ is a Cohen-Macaulay $R/x_1^dR$-module with $\dim_{R/x_1^dR} C = d-1$.

(2) $\mathrm{e}_I^0(C) = d{\cdot}\ell_R(R/I)$.

(3) $C/IC \cong (R/I)^d$.

(4) Therefore $C$ is an Ulrich $R/x_1^dR$-module with respect to $I$.

*Proof.* We may assume that $d \ge 2$. Because $\det \partial = x_1^d$, we have $x_1^d \in \mathrm{Ann}_R(C)$. Hence $C$ is a Cohen-Macaulay $(R/x_q^dR)$-module with $\dim_{R/x_q^dR} C = d-1$. We set $\mathfrak{q} = (x_2, x_3, \cdots, x_d)$ then since $I^d(R/x_1^dR) = \mathfrak{q}I^{d-1}(R/x_1^dR)$ holds true, $\mathfrak{q}$ forms a minimal reduction of $I$ in $R/x_1^dR$ so that $x_2, x_3, \cdots, x_d$ be a parameter ideal for $C$. Tensoring $R/\mathfrak{q}$ to the exact sequence ($\flat$), we get the exact sequence

$$0 \to \overline{R}^d \overset{\overline{\partial}}{\to} \overline{R}^d \to C/\mathfrak{q}C \to 0$$

of $R$-modules because $x_2, x_3, \cdots, x_d$ forms a regular sequence on $C$ where $\overline{\partial} = R/I \otimes \partial$ and $\overline{R} = R/\mathfrak{q}$. Then we have $C/\mathfrak{q}C \cong (R/I)^d$ and $IC = \mathfrak{q}C$. Thus we get $\mathrm{e}_I^0(C) = \ell_R(C/\mathfrak{q}C) = d{\cdot}\ell_R(R/I)$ and $C/IC \cong (R/I)^d$. Consequently, our required conditions are satisfied. □

**Remark 4.3.** Assume that the equality $\chi_1(Q;M) = \mathrm{hdeg}_Q(M) - \mathrm{e}_Q^0(M)$ holds true for a parameter ideal $Q$ for $M$. Then, because $\mathrm{e}_Q^1(M) = -\mathrm{T}_Q^1(M)$ by [GO, Theorem 1.3], the equality $\mathrm{g}_s(Q;M) = \mathrm{hdeg}_Q(M) - \mathrm{e}_Q^0(M) - \mathrm{T}_Q^1(M)$ is satisfied. However, the converse does not hold true in general.

We close this paper with the following example of parameter ideals $Q$ such that $\mathrm{g}_s(Q;A) = \mathrm{hdeg}_Q(A) - \mathrm{e}_Q^0(A) - \mathrm{T}_Q^1(M)$ but $\chi_1(Q;M) < \mathrm{hdeg}_Q(M) - \mathrm{e}_Q^0(M)$.

**Example 4.4.** Let $\ell \ge 2$ and $m \ge 1$ be integers. Let

$$S = k[[X_i, Y_i, Z_j \mid 1 \le i \le \ell, 1 \le j \le m]]$$

be the formal power series ring with $2\ell + m$ indeterminates over an infinite field $k$. Let

$$A = S/(X_1, X_2, \dots, X_\ell) \cap (Y_1, Y_2, \dots, Y_\ell),$$
$$\mathfrak{m} = (x_i, y_i, z_j \mid 1 \le i \le \ell, 1 \le j \le m)A, \text{ and}$$
$$Q = (x_i - y_i \mid 1 \le i \le \ell)A + (z_j \mid 1 \le j \le m)A,$$

where $x_i$, $y_i$, and $z_j$ denote the images of $X_i$, $Y_i$, and $Z_j$ in $A$ respectively. Then $\mathfrak{m}^2 = Q\mathfrak{m}$, whence $Q$ is a reduction of $\mathfrak{m}$. We furthermore have the following:

(1) $A$ is an unmixed local ring with $\dim A = \ell + m$, $\mathrm{depth} A = m+1$, and $\mathrm{H}_{\mathfrak{m}}^{m+1}(A)$ is not finitely generated.

(2) $\ell_A(A/Q) = \ell + 1$, $\mathrm{e}_Q^0(A) = 2$, $\mathrm{e}_Q^1(A) = -1$, and hence $\chi_1(Q;A) = \ell - 1$ and $\mathrm{g}_s(Q;A) = \ell - 2$.

(3) $\mathrm{hdeg}_Q(A) = 2 + \binom{\ell+m-1}{m+1}$ and $\mathrm{T}^1_Q(A) = \binom{\ell+m-2}{m}$.

(4) Hence $\mathrm{g}_s(Q;A) = \mathrm{hdeg}_Q(A) - \mathrm{e}^0_Q(A) - \mathrm{T}^1_Q(A)$, if $\ell = 2, 3$, but $\chi_1(Q;A) < \mathrm{hdeg}_Q(A) - \mathrm{e}^0_Q(A)$ if $\ell = 3$.

*Proof.* Set $\mathfrak{a}_1 = (X_i \mid 1 \le i \le \ell)$ and $\mathfrak{a}_2 = (Y_i \mid 1 \le i \le \ell)$ and consider the exact sequence

$$0 \to A \to S/\mathfrak{a}_1 \times S/\mathfrak{a}_2 \to S/[\mathfrak{a}_1 + \mathfrak{a}_2] \to 0$$

of $S$-modules. Then because

$$S/\mathfrak{a}_1 \cong k[[Y_i, Z_j \mid 1 \le i \le \ell, 1 \le j \le m]],$$
$$S/\mathfrak{a}_2 \cong k[[X_i, Z_j \mid 1 \le i \le \ell, 1 \le j \le m]], \quad \text{and}$$
$$S/[\mathfrak{a}_1 + \mathfrak{a}_2] \cong k[[Z_j \mid 1 \le j \le m]],$$

we get $\dim A = \ell + m$, $\mathrm{H}^{m+1}_{\mathfrak{m}}(A) \cong \mathrm{H}^m_{\mathfrak{m}}(S/[\mathfrak{a}_1 + \mathfrak{a}_2])$, and $\mathrm{H}^j_{\mathfrak{m}}(A) = 0$ for all $j \ne m+1$, $\ell + m$. Hence we have

$$\mathrm{hdeg}_Q(A_{m+1}) = \mathrm{hdeg}_Q(S/[\mathfrak{a}_1 + \mathfrak{a}_2]) = \mathrm{e}^0_Q(S/[\mathfrak{a}_1 + \mathfrak{a}_2]) = \mathrm{e}^0_{\mathfrak{m}}(S/[\mathfrak{a}_1 + \mathfrak{a}_2]) = 1$$

and $\mathrm{hdeg}_Q(A_j) = 0$ for all $0 \le j \le \ell + m - 1$ such that $j \ne m + 1$. Therefore, because $\mathrm{e}^0_Q(A) = \mathrm{e}^0_{\mathfrak{m}}(A) = 2$, we get

$$\mathrm{hdeg}_Q(A) = \mathrm{e}^0_Q(A) + \sum_{j=0}^{\ell+m-1} \binom{\ell+m-1}{j} \mathrm{hdeg}_Q(A_j) = 2 + \binom{\ell+m-1}{m+1}$$

and

$$\mathrm{T}^1_Q(A) = \sum_{j=1}^{\ell+m-1} \binom{\ell+m-2}{j-1} \mathrm{hdeg}_Q(A_j) = \binom{\ell+m-2}{m}.$$

Let $S' = k[[X_i, Y_i \mid 1 \le i \le \ell]]$ be the formal power series ring and we set

$$B = S'/[\mathfrak{a}_1 S' \cap \mathfrak{a}_2 S']$$

and $Q_0 = (x_i - y_i \mid 1 \le i \le \ell)B$. Then we have $A = B[[Z_j \mid 1 \le j \le m]]$ and $Q = Q_0 A + (z_j \mid 1 \le j \le m)A$. Because $z_1, z_2, \cdots, z_m$ forms a superficial sequence for $A$ with respect to $Q$ (recall that $\mathrm{gr}_Q(A) = \mathrm{gr}_{Q_0}(B)[W'_1, W'_2, \cdots, W'_m]$ forms the polynomial ring, where $W_j$'s are the initial forms of $z_j$'s) and $B$ is a Buchsbaum ring with $\mathrm{H}^1_{\mathfrak{n}}(B) \cong k$ and $\mathrm{H}^i_{\mathfrak{n}}(B) = (0)$ for all $i \ne 1, \ell$, we have $\mathrm{e}^1_Q(A) = \mathrm{e}^1_{Q_0}(B) = -1$ ([SV, Proposition 2.7]). □

## References

[AB] M. Auslander and D. Buchsbaum, *Codimension and Multiplicity*, Ann. Math. **68** (1958), 625–657.

[GHV] S. Goto, J.-Y. Hong and W. V. Vasconcelos, *The homology of parameter ideals*, J. Algebra **368** (2012), 271–299.

[GhGHOPV] L. Ghezzi, S. Goto, J. Hong, K. Ozeki, T. T. Phuong, and W. V. Vasconcelos, *The Chern and Euler coefficients of modules*, arXiv:1109.5628 [math.AC] (preprint).

[GNi] S. Goto and K. Nishida, *Hilbert coefficients and Buchsbaumness of associated graded rings*, J. Pure and Appl. Algebra **181** (2003), 61–74.

[GO] S. Goto and K. Ozeki, *The first Euler characteristics versus the homological degrees*, arXiv:1404.2455 [math.AC] (preprint).

[GOTWY] S. Goto, K. Ozeki, R. Takahashi, K. Watanabe, and K. Yoshida, *Ulrich ideals and modules*, Math. Proc. Camb. Phil. Soc. **156** (2014), 137–166.

[H] C. Huneke, *On the symmetric and Rees algebra of an ideal generated by a d-sequence*, J. Algebra **62** (1980), 268–275.

[MSV] M. Mandal, B. Singh, and J. Verma, *On some conjectures about the Chern numbers of filtrations*, J. Algebra **325** (2011), 147–162.

[Mc] L. McCune, *Hilbert Coefficients of Parameter Ideals*, arXiv:1111.4186 [math.AC] (preprint).

[N] M. Nagata, Local Rings, Interscience, 1962.

[O] A. Ooishi, *$\delta$-genera and sectional genera of commutative rings*, Hiroshima Math. J. **17**, 1987, 361–372.

[Se] J. P. Serre, *Algébre locale, Multiplicités,* Lecture Notes in Mathematics **11**, Springer, Berlin, 1965.

[STC] P. Schenzel, N. V. Trung and N. T. Cuong, *Verallgemeinerte Cohen-Macaulay-Moduln*, Math. Nachr. **85** (1978), 57–73.

[SV] J. Stückrad, W. Vogel, Buchsbaum Rings and Applications, Springer-Verlag, Berlin, Heidelberg, New York, 1986.

[V1] W. V. Vasconcelos, *The homological degree of a module*, Trans. Amer. Math. Soc. **350** (1998), 1167–1179.

[V2] W. V. Vasconcelos, *Cohomological degrees of graded modules* in "Six lectures on Commutative Algebra", Progress in Mathematics **166**, 345–392, Birkhäuser Verlag, Basel · Boston · Berlin.

Department of Mathematics, School of Science and Technology, Meiji University, 1-1-1 Higashi-mita, Tama-ku, Kawasaki 214-8571, Japan

*E-mail address*: goto@math.meiji.ac.jp

Department of Mathematical Science, Faculty of Science, Yamaguchi University, 1677-1 Yoshida, Yamaguchi 753-8512, Japan

*E-mail address*: ozeki@yamaguchi-u.ac.jp